\documentclass{article}
\usepackage[latin1]{inputenc}
\usepackage{amsmath}
\usepackage{amssymb}
\usepackage{amsthm}
\usepackage{mathrsfs}
\allowdisplaybreaks

\begin{document}

\theoremstyle{plain} \newtheorem{Thm}{Theorem}
\newtheorem{prop}{Proposition}
\newtheorem{lem}{Lemma}
\newtheorem{step}{Step}
\newtheorem{de}{Definition}
\newtheorem{obs}{Observation}
\newtheorem{cor}{Corollary}

\def\bfs{\bfseries}

\newcommand{\kor}{\mathring }
\newcommand{\Korz}{\mathring{\mathbb{Z}} }
\newcommand{\korx}{\mathring{x} }
\newcommand{\res}{\restriction }
\newcommand{\kalapz }{\hat{\mathbb{Z}} }
\theoremstyle{remark} \newtheorem*{pf}{Proof}
\renewcommand\theenumi{(\alph{enumi})}
\renewcommand\labelenumi{\theenumi}
\renewcommand{\qedsymbol}{}
\renewcommand{\qedsymbol}{\ensuremath{\blacksquare}}
\title{The Cayley isomorphism property for groups of order $p^3q$}
\author{Gábor Somlai\\
Department of Algebra and Number Theory \\
E\"otv\"os Loránd University \\ 1117 Budapest, Pázmány Péter sétány 1/C , Hungary\\
email: gsomlai@cs.elte.hu\thanks{Research supported by the Hungarian
  Scientific Fund (OTKA), grant no. K84233}}
\date{}
\maketitle

\begin{abstract}
For every prime $p > 3$ and for every prime $q>p^3$ we  prove that $\mathbb{Z}_q \times \mathbb{Z}_p^3$  is a DCI-group.

\end{abstract}

\section{Introduction}\label{intro}

Let $G$ be a finite group and $S$ a subset of $G$. The Cayley graph
$Cay(G,S)$ is defined by having the vertex set $G$ and $g$ is adjacent
to $h$ if and only if $g h^{-1} \in S$. The set $S$ is called the
connection set of the Cayley graph $Cay(G,S)$. A Cayley graph
$Cay(G,S)$ is undirected if and only if $S=S^{-1}$, where $S^{-1} =
\left\{ \, s^{-1} \in G \mid s \in S \, \right\} $.
Every right multiplication via elements of $G$ is an automorphism of
$Cay(G,S)$, so the automorphism group of every Cayley graph on $G$
contains a regular subgroup isomorphic to $G$. Moreover, this property
characterises the Cayley graphs of $G$.

It is clear that automorphism $\mu$ of the group $G$ induces an isomorphism between $Cay(G,S)$ and $Cay(G,S^{\mu } )$.
Such an isomorphism is called a Cayley isomorphism. A Cayley graph $Cay(G,S)$ is said to be a
CI-graph if, for each $T \subset G$, the Cayley graphs $Cay(G,S)$ and
$Cay(G,T)$ are isomorphic if and only if there is an automorphism
$\mu$ of $G$ such that $S^{\mu } =T$. Furthermore, a group $G$
is called a DCI-group if every Cayley graph of $G$ is a CI-graph and
it is called a CI-group if every undirected Cayley graph of $G$ is a CI-graph.

The problem of investigating the isomorphism problem of Cayley graphs started with Ádám's conjecture \cite{adam}, which states that every circulant graph if a CI-graph. Using our terminology, it was conjectured that every cyclic group is a DCI-group. This conjecture was first disproved by Elspas and Turner \cite{ElspasTurner} for directed Cayley graphs of $\mathbb{Z}_8$ and for undirected graphs of Cayley graphs of $\mathbb{Z}_{16}$.

By investigating the spectrum of circulant graph Elspas and Turner \cite{ElspasTurner}, and independently Djokovi\'{c} \cite{Djokovic} proved that every cyclic group of order $p$ is a CI-group if $p$ is a prime. Also a lot of research was devoted to the investigation of circulant graphs. One of the most important results for our investigation is that $\mathbb{Z}_{pq}$ is a DCI-group for every pair of primes $p<q$. This result was first proved by Alspach and Parsons \cite{AlspachParsons} and later by Pöschel and Klin \cite{KlinPochel4} using Schur rings, and by Godsil \cite{godsil}. Finally, Muzychuk \cite{Muzysquarefree,Muzycyclic1} proved that a cyclic group $\mathbb{Z}_n$ is a DCI-group if and only if $n=k$ or $n=2k$, where $k$ is square-free. Furthermore, $\mathbb{Z}_n$ is a CI-group if and only if $n$ is as above or $n=8,9,18$.

It is easy to see that every subgroup of a (D)CI-group is also a (D)CI-group so it is natural to investigate $p$-groups which are the Sylow $p$-subgroups of a finite group. Babai and Frankl \cite{babaifrankl} proved that if $G$ is a $p$-group, which is a CI-group, then $G$ can only be elementary abelian $p$-group, the quaternion group of order $8$ or one of a few cyclic groups $\mathbb{Z}_4$, $\mathbb{Z}_8$, $\mathbb{Z}_9$ or $\mathbb{Z}_{27}$.
Muzychuk's result about cyclic groups shows that $\mathbb{Z}_{27}$ is not a CI-group and $\mathbb{Z}_8$ is not a DCI-group. They also asked whether every elementary abelian $p$-group is a CI-group.

The cyclic group of order $p$, which is a CI-group, can also be considered as an elementary abelian $p$-group of rank $1$. The best general result  was given by Hirasaka and Muzychuk \cite{HirasakaMuzy} who proved that $\mathbb{Z}_p^4$ is a CI-group for every prime $p$. For our investigation the following weaker results are also important. Dobson \cite{dobson} proved that $\mathbb{Z}_p^3$ is a CI-group for every prime $p$ and Alspach and Nowitz shoved \cite{alspachnowitz} that $\mathbb{Z}_p^3$ is a CI-group with respect to Cayles color digraphs. However Muzychuk \cite{Muzychukelemuabel} showed that an elementary abelian $p$-group of $2p-1 + \binom{2p-1}{p}$ rank is not a CI-group.

Severe restriction on the structure of CI-groups was given by Li and Praeger and then a more precise list of candidates for CI-groups was given by Li, Lu and Pálfy \cite{liluppp}.

New family of CI-groups was found by Kovács and Muzychuk \cite{KovacsMuzy}, that is, $\mathbb{Z}_{p^2} \times \mathbb{Z}_q$ is a CI-group for every prime $p$ and $q$. It was also conjectured in \cite{KovacsMuzy}, that the direct product of CI-groups of coprime order is a CI-group.

\begin{Thm}\label{fotetel}
For every prime $p$ and every prime $q>p^3$ the group $\mathbb{Z}_{p^3} \times \mathbb{Z}_q$ is a DCI-group.
\end{Thm}

Our paper is organized as follows. In Section \ref{technicaldetails} we introduce the notation that will be used throughout this paper. In Section \ref{basicideas} we collect important ideas that we will use in the proof of Theorem \ref{fotetel}. Finally, Section \ref{mainresult} contains the proof of Theorem \ref{fotetel}.

\section{Technical details}\label{technicaldetails}
In this section we introduce some notation.
Let $G$ be a group. We use $H \le G$ to denote that $H$ is a subgroup
of $G$ and by $N_{G} (H)$ and $C_{G} (H)$ we denote the normalizer and
the centralizer of $H$ in $G$, respectively. The center of a group $G$
will be denoted by $Z(G)$.

Let us assume that the group $H$ acts on the set $\Omega$ and let $G$ be an arbitrary group. Then by $G \wr_{\Omega} H$ we denote the wreath
product of $G$ and $H$. Every element $g \in G \wr_{\Omega} H$ can be
uniquely written as $hk$, where $k \in K = G^{\Omega}$ and $h \in H$. The group $K = G^{\Omega}$ is
called the base group of $G \wr_{\Omega} H$ and the elements of $K$ can
be treated as functions from $\Omega$ to $G$. If $g \in G \wr_{\Omega} H$ and $g=hk$
we denote $k$ by $(g)_b$. In order to simplify the notation $\Omega$ will be omitted if it is clear from
the definition of $H$ and we will write $G \wr H$.

The symmetric group on the set $\Omega$ will be denoted by $Sym(\Omega)$. Let $G$ be a permutation group on the
set $\Omega $. For a $G$-invariant partition $\mathcal{B}$ of the set
$\Omega$ we use $G^{\mathcal{B}}$ to denote the permutation group on $\mathcal{B}$ induced by the action of $G$ and
similarly, for every $g \in G$ we denote by $g^{\mathcal{B}}$ the
action of $g$ on the partition $\mathcal{B}$.

For a group $G$, let $\hat{G}$ denote the subgroup of the
symmetric group $Sym(G)$ formed by the elements of $G$ acting by
right multiplication on $G$. For every Cayley graph $\Gamma =Cay(G,S)$ the
subgroup $\hat{G}$ of $Sym(G)$ is contained in $Aut(\Gamma)$.

\begin{de}\label{zart}
Let $G \le Sym(\Omega)$ be a permutation group. Let
\begin{equation*}
\begin{split}
G^{(2)} = \left\{ \pi \in Sym(\Omega) \bigg| \begin{matrix}
\forall a,b \in
  \Omega \mbox{  } \exists g_{a,b} \in G \mbox{ with}  &\pi
  (a)=g_{a,b} (a)
 \mbox{ and } \\ &\pi (b)=g_{a,b} (b) \end{matrix} \right\} \mbox{.}
\end{split}
\end{equation*}
We say that $G^{(2)}$ is the $2$-closure of the permutation group $G$.
\end{de}

\begin{lem}\label{2zart}
Let $\Gamma$ be a graph. If $G \le Aut(\Gamma )$, then
$G^{(2)} \le Aut(\Gamma )$.
\end{lem}

\section{Basic ideas}\label{basicideas}
In this section we collect some results and some important ideas that
we will use in the proof of Theorem \ref{fotetel}.

We begin with a fundamental lemma that we will use all along this paper.
\begin{lem}[Babai \cite{babai}]\label{babai}
$Cay(G,S)$ is a CI-graph if and only if for every regular subgroup
$\mathring{G}$ of $Aut(Cay(G,S) )$ isomorphic to $G$ there is a
$\mu \in Aut(Cay(G,S))$ such that ${\mathring{G}}^{\mu} = \hat{G}$.
\end{lem}

We introduce the following definition.
\begin{de}
\begin{enumerate}
\item
We say that a Cayley graph $Cay(G,S)$ is a CI$^{(2)}$-graph if and only if for every regular subgroup $\mathring{G}$ of
$Aut(Cay(G,S) )$ isomorphic to $G$ there is a $\sigma \in \langle \mathring{G}, \hat{G} \rangle ^{(2)}$
such that ${\mathring{G}}^{\sigma} = \hat{G}$.
\item A group $G$ is called a DCI$^{(2)}$-group if for every $S \subset G$ the Cayley graph $Cay(G,S)$ is a CI$^{(2)}$-graph.
\end{enumerate}
\end{de}

\begin{de}
Let $\Gamma$ be an arbitrary graph and $A, B \subset V(
\Gamma)$ such that $A \cap B = \emptyset$. We write $A \sim B$ if one of the following four possibilities holds:
\begin{enumerate}
\item For every $a \in A$ and $b \in B$ there is an edge from $a$ to
      $b$ but there is no edge from $b$ to $a$.
\item For every $a \in A$ and $b \in B$ there is an edge from $b$ to
      $a$ but there is no edge from $a$ to $b$.
\item For every $a \in A$ and $b \in B$ the vertices $a$ and $b$ are
      connected with an undirected edge.
\item There is no edge between $A$ and $B$.
\end{enumerate}
We also write $A \nsim B$ if none of the previous four possibilities holds.
\end{de}

\begin{lem}\label{koszoru}
Let $A$, $B$ be two disjoint subsets of cardinality $p$ of a
graph. We write $A \cup B = \mathbb{Z}_p \cup \mathbb{Z}_p $. Let
us assume that $\hat{\mathbb{Z}}_p$ acts naturally on $A \cup B$
and for a generator $\mathring{g} $ of the cyclic group
$\mathring{Z}_p$ the action of $\mathring{a}$ is defined by
$(a_1,a_2)\mathring{g}=(a_1 +b, a_2 +c)$ for some $b, c \in
\mathbb{Z}_p$.
\begin{enumerate}
\item\label{kosza} If $b=c$, then the action of $\hat{\mathbb{Z}}_p$ and
      $\mathring{Z}_p$ on $A \cup B$ are the same.
\item\label{koszb} If $A \nsim B$, then $b=c$.
\item\label{koszc} If $A \sim B$, then every $\pi \in Sym(A \cup B)$
      which fixes $A$ and $B$ setwise is an automorphism of the graph
      defined on $A \cup B$ if $\pi \restriction A \in Aut(A)$ and
      $\pi \restriction B \in Aut(B)$.
\end{enumerate}
\end{lem}
\begin{pf}
These statements are obvious.
\qed
\end{pf}

\begin{lem}\label{centrum}
Let us assume that $H$ is a regular abelian subgroup of $Sym(p^n)$ and let $P \ge H$ be a Sylow $p$-subgroup of $Sym(p^n)$.
Then $H$ contains $Z(P)$.
\end{lem}
\begin{pf}
It is well known that the center of $P$ is a cyclic $p$-group.
Let $z$ be a generator of $Z(P)$. Then $\langle H, z \rangle$ is a transitive abelian group.
Hence $\langle H, z \rangle$ is regular. Since $H$ is also regular, we have that
$z$ has to be in $H$.
\qed
\end{pf}
\section{Main result}\label{mainresult}

In this section we will prove that $\mathbb{Z}_p^3  \times
\mathbb{Z}_q$ is a DCI-group if $q>p^3$ and $p>3$.

Our technique is based on Lemma \ref{babai} so we fix a
Cayley graph $\Gamma =Cay(\mathbb{Z}_p^3 \times \mathbb{Z}_q,S)$. Let
$A=Aut(\Gamma) $ and $ \mathring{G}=\mathring{\mathbb{Z} }_p^3 \times
\mathring{\mathbb{Z} }_q $ be a regular subgroup of $ A $ isomorphic
to $\mathbb{Z}_p^3 \times \mathbb{Z}_q$. In order to prove Theorem
\ref{fotetel} we  have to find an $\alpha \in A$ such that
$\mathring{G}^{\alpha} = \hat{G} =\hat{\mathbb{Z}_p^3} \times \hat{\mathbb{Z} }_q $,
what we will achieve in three steps.

\subsection{Step 1}\label{step1}

We may assume $\hat{\mathbb{Z}}_q$ and
$\mathring{\mathbb{Z}}_q$
lie in the same Sylow $q$-subgroup $Q$ of $Sym(p^3q)$. Then both
$\mathring{\mathbb{Z}}_p^3$ and $\hat{\mathbb{Z}}_p^3$ are
subgroups of $N_{Sym(p^3q) } (Q) \cap A$ so we may assume that
$\mathring{\mathbb{Z}_p^3}$ and $\hat{\mathbb{Z}_p^3}$ lie in
the same Sylow $p$-subgroup of $N_{Sym(p^3q) } (Q) \cap A$ which
is contained in a Sylow $p$-subgroup $P$ of $A$.

The Sylow $q$-subgroup $Q$ gives a partition $\mathcal{B}=
\left\{ B_1, B_2, \ldots ,B_{p^3}
\right\}$ of the vertices of $\Gamma$, where $\left| B_i \right|
=q$ for every $i=1, \ldots ,p^3$.
It is easy to see that $\mathcal{B}$ is invariant under the action of
$\hat{\mathbb{Z}}_p^3$ and $\mathring{\mathbb{Z}}_p^3$ and hence
$\langle \hat{G} , \mathring{G} \rangle \le Sym(q) \wr Sym(p^3)$.
Moreover, both  $\mathring{G}$ and $\hat{G}$ are regular so $\mathring{\mathbb{Z}_p^3}$ and $\hat{\mathbb{Z}_p^3}$ induce regular
action on $\mathcal{B}$ which we denote by $H_1$ and $H_2$, respectively. The assumption that $\mathring{\mathbb{Z}_p^3}$ and $\hat{\mathbb{Z}_p^3}$ lie in the same Sylow $p$-subgroup of $A$
 implies that $H_1$ and $H_2$ are in the same Sylow $p$-subgroup of $Sym(p^3)$, what
 we denote by $P_1$.

\subsection{Step 2}\label{step2}
Let us assume that $\hat{\mathbb{Z} }_q \ne \mathring{\mathbb{Z}}_q$
which is generated by $p^3$ disjoint $q$-cycles.
We intend to find an element $\alpha \in A$ such that
$\mathring{\mathbb{Z}}_q^{\alpha}= \hat{\mathbb{Z}}_q$.

We define a graph $\Gamma _0$ on $\mathcal{B}$ such that $B_i$ is
connected to $B_j$ if and only if $B_i \nsim B_j$. This is an undirected
graph with vertex set $\mathcal{B}$ and both $\mathring{\mathbb{Z}}_p^3$ and
$\hat{\mathbb{Z} }_p^3$ are regular subgroups of $Aut(\Gamma_0)$.
It follows that $\Gamma_0$ is a Cayley graph of $\mathbb{Z}_p^3$.

\begin{de}
\begin{enumerate}
\item For a pair $(B_i, B_j) \in \mathcal{B}^{2}$ we write $B_i \equiv B_j$ if either there exists a path
$C_1, C_2, \ldots ,C_n$ in $\Gamma_0$ such that $C_1 =B_1$, $C_n =B_2$ or $i=j$.
\item For a pair $(B_i, B_j) \in \mathcal{B}^{2}$ we write $B_i \not\equiv B_j$ if $B_i \equiv B_j$ does not hold.
\item If both $H$ and $K$ are subsets of the vertices of $\Gamma_0$ such that $H \cap K = \emptyset$ and
for every $B_i \in H$, $B_j \in K$ we have $B_i \not\equiv B_j$, then we write $H \not\equiv K$.
\end{enumerate}
\end{de}

\begin{obs}
\begin{enumerate}
\item
The relation $\equiv$ defines an equivalence relation on $\mathcal{B}$.
The connected components of $\Gamma_0$ will be called equivalence classes.
\item Since $H_1$ acts transitively on $\mathcal{B}$
we have that the size of the equivalence classes defined by the relation $\equiv$ divides $p^3$.
\end{enumerate}
\end{obs}

We can also define a colored graph $\Gamma_1$ on $\mathcal{B}$
by coloring the edges of the complete directed graph on $p^3$
points. $B_i$ is connected to $B_j$ with the same color as $B_i'$ is connected to $B_j'$ in $\Gamma_1$ if and only if there exists a graph isomorphism $\phi$ from $B_i \cup B_j$ to $B_i' \cup B_j'$ such that $\phi(B_i) = B_i'$ and $\phi(B_j) = B_j'$. The graph $\Gamma_1$
is a colored Cayley graph of the elementary abelian $p$-group
$\mathbb{Z}_p^3$. Moreover, both $H_1$ and $H_2$ act regularly
on $\Gamma_1$.

We prove the following two lemmas what we will use several times
in this step.

\begin{lem}\label{aut}
Let us assume that $C_1', C_2', \ldots C_k'$ are the equivalence classes defined in $V(\Gamma_0)$ and let $C_i = \cup C_i' \subset V(\Gamma)$ for every $i=1, \ldots , k$ . Let $\alpha$ be a permutation on the vertex set $V(\Gamma)$ such that for every $1 \le i \le k$ the restriction $\alpha \restriction C_i=\eta_i \restriction C_i$ for some $\eta_i \in Aut(\Gamma)$ and $\alpha^{  V(\Gamma_0)   }$ is an automorphism of $\Gamma_0$. Then $\alpha$ is an automorphism of $\Gamma$.
\end{lem}
\begin{pf}
Let $x$ and $y$ be points in $V(\Gamma)$. We have to prove that $x$ is connected to $y$ if and only if $\alpha(x)$ is connected to $\alpha(y)$.
This holds if $x$ and $y$ are in the same $C_i$ for some $1\le i \le k$ since $\alpha \restriction C_i$ is defined by an automorphism of $\Gamma$ on $C_i$.
If $x \in B_m $ and $y \in B_n$, where $B_m \sim B_n$ and $x$ is connected to $y$, then every element of $B_m$ is connected to every element of $B_n$. Since $\alpha^{ V(\Gamma_0) } \in Aut(\Gamma_0)$ the same holds for $\alpha(B_m)$ and $\alpha(B_n)$ and hence $\alpha(x)$ is connected to $\alpha(y)$. Similar argument shows that if $x \in B_m $ and $y \in B_n$, where $B_m \sim B_n$ and $x$ is not connected to $y$, then $\alpha(x)$ is not connected to $\alpha(y)$.
\qed
\end{pf}

\begin{lem}\label{szorzas}
\begin{enumerate}
\item
Let $A$ and $B$ be two disjoint subsets of cardinality $q$ of $V(\Gamma)$.
We write $A= \{(a,x) \mbox{ }| \mbox{ }x \in \mathbb{Z}_q )\}$ and $B =\{(b,x) \mbox{ }| \mbox{ }x \in \mathbb{Z}_q \}$. Let us assume that
$\hat{g}$ and $\mathring{g}$ are automorphisms of the graph $\Gamma$ with $\hat{g}(a,x)=\mathring{g}(a,x)=(a,x+1)$,
$\hat{g}(b,x)=(b,x+1)$ and $\mathring{g}(b,x)=(b,x+d)$ for some $d \in \mathbb{Z}_q$ for all $x \in \mathbb{Z}_q$.
 Furthermore, let us assume that $\hat{w}$ and $\mathring{w}$ are automorphisms of
the graph $\Gamma$ with $\hat{w}(A)=\mathring{w} (A)=B$ and $\hat{w}$ and $\mathring{w}$ commute with $\hat{g}$ and
$\mathring{g}$, respectively. Then for $\alpha=\mathring{w} \hat{w}^{-1}$ we have $\mathring{g}^{\alpha}\restriction_B
=\hat{g} \restriction_B$.
\item
Let us assume that $C = \{ (c,x) \mid x \in \mathbb{Z}_q \}$ is a subset of $V(\Gamma)$  with $A \cap B= A \cap C=\emptyset$.
 We also assume that $\hat{g}(c,x) =(c,x+1)$ and $\mathring{g}(c,x) =(c,x+d)$ for every $x\in \mathbb{Z}_q$.
 Let us assume that $\mathring{v} \in Aut(\Gamma)$ with $\mathring{v}(A)=C$ and we also assume that $\mathring{g}$ and $\mathring{v}$ commute.
 Then for $\beta =\mathring{v} \hat{w}^{-1}$ we have $\mathring{g}^{\beta}\restriction_B
=\hat{g} \restriction_B$.

\end{enumerate}
\end{lem}
\begin{pf}
\begin{enumerate}
\item\label{szorzasa}
Let us assume that $\hat{w} (a,0)=(b,b_0)$  and $\mathring{w} (a,0)=(b,b_0')$ for some $b_0, b_0' \in \mathbb{Z}_q$.
 Using that $\hat{w}$ and
 $\hat{g}$ commute we get that $\hat{w}(a,x) = (b,b_0+x)$ for every $x \in \mathbb{Z}_q$ and similarly we have
 $\mathring{w}(a,x)=(b,b_0'+dx)$.
 Thus
\begin{equation*} \begin{split}
  \alpha\left(b,x\right)&=\alpha\left(b,b_0+(x-b_0) \right) =\mathring{w} \left(a,x-b_0 \right)=
 \left(b,b_0' +(x-b_0)d \right)\\&=\left(b,(b_0' -db_0) +dx \right) \mbox{.}  \end{split} \end{equation*}
It is easy to derive that $\alpha^{-1}(b,x) =\left( b, \frac{x - (b_0'-db_0)}{d} \right)$.
Using the previous two equations we get
\begin{equation*} \begin{split} &\alpha^{-1} \mathring{g} \alpha  \restriction_B (b,x) = \alpha^{-1} \mathring{g}
\left(b,(b_0' -db_0) +dx \right) = \alpha^{-1}\left(b,(b_0' -db_0) +dx+d\right) \\&=
\left(b, \frac{(b_0' -db_0) +d x + d -(b_0'-db_0)}{d} \right)=(b,x+1) \mbox{.}\end{split} \end{equation*}
\item\label{szorzasb}
Let us assume that $\mathring{v}(a,0) =(c,c_0)$ for some $c_0 \in \mathbb{Z}_q$.
 Then
 $\mathring{v}(a,x) =(c,c_0 +dx)$ for all $x \in \mathbb{Z}_q$
Thus
 \begin{equation*} \begin{split} \beta(b,x) &=\mathring{v} \hat{w}^{-1} \left( b,b_0+(x-b_0)\right) \\ &= \mathring{v}\left(a,x-b_0 \right)
 =\left(c,c_0+(x-b_0)d \right)\end{split} \end{equation*} and hence
$\beta^{-1} (c,x) =(b,\frac{x-c_0+b_0 d}{d})$.
Similarly to the previous case we have
\begin{equation*} \begin{split}
\beta^{-1} \mathring{g} \beta \left( b,x \right) &= \beta^{-1} \mathring{g} \left( c, c_0+(x-b_0)d \right)
= \beta^{-1} \left( c, c_0+(x-b_0)d + d \right) \\&= \left( b, \frac{c_0+(x-b_0)d + d -c_0 + b_0 d}{d}\right) =\left(b,x+1 \right) \mbox{.}
\end{split} \end{equation*}
\end{enumerate}
\qed
\end{pf}

The points of the graph $\Gamma_0$ and $\Gamma_1$ can be identified
with the elements of $\mathbb{Z}_p^3$ and we may assume that the action
of an element $r$ of the Sylow $p$-subgroup $P_1$ is the following:
\[ r(a,b,c) =(a+x, b+s_a, c+t_{a,b}) \mbox{,}\]
where $s_a$ only depends on $a$ and $t_{a,b}$ depends on $a$ and $b$.

Let $\hat{g}$ and $\mathring{g}$ denote the generator of
$\hat{\mathbb{Z}}_q$ and $\mathring{\mathbb{Z}}_q$, respectively.
We may assume that $\hat{g}\restriction B_1 = \mathring{g}\restriction B_1$.

\begin{enumerate}
\item
Let us assume first that $\Gamma_0$ is a connected graph.

Using Lemma \ref{koszoru} \ref{koszb} we get that $\hat{g} \res{B_i} = \mathring{g} \res{B_j}$
 if there exists a path in $\Gamma_0$ from $B_i$ to $B_j$.
This shows that $\hat{g}=\mathring{g}$ since $\Gamma_0$ is connected in this case.

\item
Let us assume that $\Gamma_0$ is the empty graph.

 For every $B_m \in \mathcal{B}$ there exist $\hat{r}_m$
 and $\mathring{r}_m$ such that $\hat{r}_m (B_1) =
 \mathring{r}_m (B_1)=B_m$.

 Let $\alpha$ be defined as follows
\begin{equation} \begin{split}
\alpha \restriction B_1 &= id
\\ \alpha \restriction B_m &= \mathring{r}_m \hat{r}_m^{-1}
\mbox{ } \mbox{ for } \mbox{ } 2 \le m \le p^3 \mbox{.}
\end{split}
\end{equation}
It is easy to see that $\alpha^{\mathcal{B}} = id$ so using Lemma \ref{aut} we get that  $\alpha $ is an automorphism of $\Gamma$.
Using Lemma \ref{szorzas} \ref{szorzasa} we get that $\mathring{g}^{ \alpha}=\hat{g}$.

\item
Let us assume that the size of the connected components
of $\Gamma_0$ is $p$.

Let $C_1', C_2', \ldots ,C_{p^2}'$ denote the equivalence
classes defined by the relation $\equiv$ on $\Gamma_0$ and for $1 \le m \le p^2$ let $C_m=\cup C_m'$. For $C_2, \ldots, C_{p^2}$ we choose
 an element $\hat{u}_m$ of $\hat{\mathbb{Z}}_p^3$ such that $\hat{u}_m(C_1) = C_m$.
 We may assume that $B_1 \subset C_1$.
Since $H_2$ is regular on $\Gamma_0$, for every $2 \le m \le p^2$
there exists $\mathring{u}_m$ such that $\mathring{u}_m(B_1) =
\hat{u}_m(B_1)$. For $2 \le m \le p^2$ let $\tilde{u}_m =\mathring{u}_m \hat{u}_m^{-1} $.
Now we define the following permutation:
\begin{equation*} \begin{split} \alpha_1 \restriction C_1 &= id \\
 \alpha_1 \restriction C_m &= \tilde{u}_m \mbox{ for } 2 \le m \le p^2.
  \end{split} \end{equation*}
Clearly, for $2 \le m \le p^2$ we have $\tilde{u}_m(B_j)=B_j$ for at least one $B_j \subset C_m$. Since $H_1$ and $H_2$ are in the same
Sylow $p$-subgroup of $Sym(p^3)$ the order of $\tilde{u}_m ^{\mathcal{B}}$ is a power of $p$.
We also have that $C_m$ is the union of $p$ elements of $\mathcal{B}$ for $1 \le m \le p^2$ hence
$\alpha_1^{\mathcal{B}} =id$.
We also have that
$\alpha_1 \restriction C_m$ is the restriction of an automorphism
 of the graph $\Gamma$ for $m=1, \ldots p$. Therefore by Lemma \ref{aut} $\alpha_1$ is an automorphism of the
 graph $\Gamma$.

Finally, Lemma \ref{szorzas} \ref{szorzasb} gives $\mathring{g}^{\alpha_1} =\hat{g}$.

\item
Let us assume that the size of the connected components of $\Gamma_0$
and hence the size of the equivalence classes is $p^2$.
Let $D_0', D_1', \ldots ,D_{p-1}'$ denote the equivalence classes and let
$D_m = \cup D_m'$ for $0 \le m \le p-1$.

Using Lemma \ref{centrum} we get that $H_1 \cap H_2 \ne \{ 1\}$.
Let $z$ be an element of order $p$ of $H_1 \cap H_2$ and we denote by $z_1$ and
$z_2$ the element of $\hat{\mathbb{Z}}_p^3$ and $\mathring{\mathbb{Z}}_p^3$
such that $z_1^{\mathcal{B}} = z_2^{\mathcal{B}} =z$, respectively.
Then $(z_2^{-i} z_1^{i})^{\mathcal{B}} =id$ for $i=1, \dots ,p-1$.

Let us assume first that $z_1(D_0) \ne D_0$. We may assume that
$z_1^{i}(D_0)=D_{i}$ for $i=0,1,\ldots, p-1$.
 We define $\alpha_2$ in the following way:
\begin{equation*} \begin{split}
\alpha_2 \restriction D_0&=id \\
\alpha_2 \restriction D_{i}&=z_2^{i} z_1^{-i} \mbox{ for } 1 \le i \le p-1
\mbox{.}\end{split}\end{equation*}

Since $z_1^{\mathcal{B}} = z_2^{\mathcal{B}} =z$ we have $\alpha_2^{\mathcal{B}}=id$.
Using Lemma \ref{aut} again we get that $\alpha_2
\in Aut(\Gamma)$ and Lemma \ref{szorzas} gives $\mathring{g}^{\alpha_2}=\hat{g}$.

Therefore we may assume that $z_1(D_0)=D_0$. In this case the
orbits of $z$ give a $\langle H_1, H_2 \rangle$-invariant
partition $\mathcal{E}=\{E_{a,b} \mid a,b \in \mathbb{Z}_p \}$
 of $\mathcal{B}$. Using that the elements of $\mathcal{B} =V(\Gamma_0)$
 can be identified with elements of $\mathbb{Z}_p^3$ we may assume that $E_{a,b}$
 has the following form for every pair $(a,b) \in \mathbb{Z}_p^2$:
  \[ E_{a,b}=\{ (a,b,c) \in \mathbb{Z}_p^3 \mid c
\in \mathbb{Z}_p \}. \] We may also assume that $D_a' = \cup_{b \in
\mathbb{Z}_p} E_{a,b}$ for all $a \in \mathbb{Z}_p$.

Since $H_1$ acts regularly on
$\Gamma_0$, there exists $h_1 \in H_1$ such that $h_1(E_{0,0}) = E_{0,1}$. Since $H_2$ is also regular, there exists
$h_2 \in H_2$ such that $h_2(E_{0,0}) =h_1(E_{0,0})$.
Since the order of
$h_1$ and $h_2$ are $p$ and $h_1(D_0')=h_2(D_0')=D_0'$
we have that $h_1(D_i')=h_2(D_i')=D_i'$ for $i=0, \ldots, p-1$.

We may assume that $z$, $h_1$ and $h_2$ act in the following way on
$\mathbb{Z}_p^3$.
\begin{equation*} \begin{split}  z(a,b,c) &=(a,b,c+1) \\
 h_1(a,b,c) &=(a,b+1,c )\\
 h_2(a,b,c) &=(a,b+s_a,c +t_{a,b}) \mbox{.}
\end{split}\end{equation*}
The assumption that $h_1(E_{0,0})=h_2(E_{0,0})=E_{0,1}$ gives that $s_0=1$.

We claim that $s_a=1$ for $1 \le a \le p-1$.
Since $H_2$ is regular on $\Gamma_0$ there exists $k_2 \in H_2$
 such that $k_2(0,0,0)=(a,0,0)$. Since $h_2$ and $k_2$ commute
we have that $k_2(0,i,0)=(a,s_ai,w_i)$ for some $w_i \in \mathbb{Z}_p$.
If $s_a \ne 1$, then such an element cannot be in the Sylow $p$-subgroup $P_1$.

Therefore $h_2(a,b,c)=(a,b+1,c+t_{a,b})$ for all $(a,b,c) \in \mathbb{Z}_p^3$,
where $t_{a,b} \in \mathbb{Z}_p$ only depends on $a$ and $b$.

\begin{lem}\label{koszorueab}
Let $a \ne a'$ be elements of $\mathbb{Z}_p$ and we fix two
more elements $b$ and $b'$ of $\mathbb{Z}_p$. Then either
$E_{a,b} \sim E_{a',b'}$ or $t_{a,b+n} =t_{a',b'+n}$ for
all $n \in \mathbb{Z}_p$.
\end{lem}
\begin{pf}
For all $m \in \mathbb{Z}_p$ the permutation $h_2^m h_1^{-m}$
fixes $E_{a,b}$ and $E_{a',b'}$. Moreover,
\begin{equation} \begin{split} h_2^{m} h_1^{-m}(a,b,c) &=
(a,b,c+ \sum_{i=1}^{n} t_{a,b-i}) \mbox{ and } \\
h_2^{m} h_1^{-m}(a',b',c)&=
(a',b', c+ \sum_{i=1}^{n} t_{a',b'-i})
\end{split} \end{equation}
One can see using Lemma \ref{koszoru} \ref{koszb} that if
$\sum_{i=1}^{n} t_{a,b-i} \ne \sum_{i=1}^{n}
t_{a',b'-i}$ for some $m \in \mathbb{Z}_p$, then $E_{a,b} \sim
E_{a',b'}$. If $\sum_{i=1}^{n} t_{a,b-i} = \sum_{i=1}^{n}
t_{a',b'-i}$ for all $m \in \mathbb{Z}_p$, then $t_{a,b+n} =
t_{a',b'+n}$ for $n \in \mathbb{Z}_p$.
\qed
\end{pf}
For each $a \in \mathbb{Z}_p$ we define the following function from $\mathbb{Z}_p$ to $\mathbb{Z}_p$:
\[ t'_a(b) := t'_{a,b} \mbox{.}\]

\begin{lem}\label{lem8}
Let us assume that $t_a(b+n)= t_a'(b'+n)$ for all $n \in
\mathbb{Z}_p$ and we denote by $k_2$ the unique element of
$H_2$ which maps $(a,b,0)$ to $(a',b',0)$.
Then $k_2(a,b+d,e) =(a',b'+d,e)$ for all $d, e \in \mathbb{Z}_p$.
\end{lem}
\begin{pf}
Since $k_2$ and $z$ commute we have $k_2(a,b,m)=(a',b',m)$ for all
$m \in \mathbb{Z}_p$. We also have that $k_2$ and $h_2$ commute which
gives $k_2(a,b+d,e)=(a',b'+d,e)$ for all $d, e \in \mathbb{Z}_p$.
\qed
\end{pf}
\begin{cor}\label{cor1}
If the conditions of Lemma \ref{lem8} hold and $k_1$ is the unique
element of $H_1$ such that $k_1(a,b,0)=(a',b',0)$, then
$k_1\restriction_{E_{a,b}} =k_2\restriction_{E_{a,b}}$.

\end{cor}
We define an equivalence relation on the set $\{ D_0',D_1', \ldots, D_{p-1}'\}$.
We write $D_a' \doteq D_{a'}'$ if and only if there exist $b$ and $b'$
in $\mathbb{Z}_p$ such that $t_{a,b+n} =t_{a',b'+n}$ for all $n \in
\mathbb{Z}_p$.

Now we can choose a point $(a,b_a,0)$ in every $D_a'$ such that
if $D_a \doteq D_{a'}$, then $t_{a,b_a+n} = t_{a',b_{a'}+n}$ for all
$n \in \mathbb{Z}_p$. For every $1 \le a \le p-1$ there exist
$\hat{v}_a \in \hat{\mathbb{Z}}_p^3$ and $\mathring{v}_a \in
\mathring{\mathbb{Z}}_p^3$ such that $\hat{v}^{\mathcal{B}}_a(0,b_0,0)
= \mathring{v}^{\mathcal{B}}_a(0,b_0,0) =(a,b_a,0)$ since    both
$H_1$ and $H_2$ are regular.

Now we can define the following permutation:
\begin{equation*} \begin{split}
\alpha_3 \restriction_{D_0}&= id \\
\alpha_3 \restriction_{D_a}&= \mathring{v}_a \hat{v}_a^{-1}
\mbox{ } \mbox{ for } 1 \le a \le p-1 \mbox{.}\\
\end{split} \end{equation*}

\begin{lem}\label{lem9}
$\alpha_3$ is an automorphism of $\Gamma$.
\end{lem}
\begin{pf}
We prove that $\alpha_3^{\mathcal{B} }$ is an automorphism of the graph $\Gamma_1$.
If $B_i \cup B_j$ is contained in $D_a'$ for some $a \in \mathbb{Z}_p$,
then $\alpha_3$ is defined by the restriction of an automorphism of $\Gamma$.
Therefore we only have to investigate those pairs $B_i, B_j$ of points which are
not in the same set $D_a'$ for any $a\in \mathbb{Z}_p$.

Let us assume that $B_i \in E_{a,b}$ and $B_j \in E_{a',b'}$.
By the definition of $\alpha_3 $, for every $c \in \mathbb{Z}_p$
at least one ${E_{c,d}}$ is fixed by $\alpha_3^{\mathcal{B} }$. Therefore $\alpha_3^{\mathcal{B} }$ fixes every
set $E_{c,d}$ since the order of $\alpha_3^{\mathcal{B} } \restriction_{D_c'}$
 is a power of $p$ for every $c\in \mathbb{Z}_p$.

 Let us assume first that $D_a \nsim D_a'$. Lemma
 \ref{koszorueab} gives that $B_i$ is connected to $B_j$
 if and only if $\alpha_3'(B_i)$ is connected to $\alpha_3'(B_j)$
 since $E_{a,b} \sim E_{a',b'}$.

 Let us now assume that $D_a' \sim D_{a'}'$. We denote by the pair
 $(\mathring{v}_a \hat{v}_a^{-1},\mathring{v}_a \hat{v}_a^{-1} )$
 the restriction of the action of $\alpha_3$ to ${D_a' \cup D_{a'}'}$.
 Since $\mathring{v}_a$ and $\hat{v}_a^{-1}$ are automorphisms of $\Gamma$
 the pair $((\mathring{v}_a \hat{v}_a^{-1})^{\mathcal{B}},
 (\mathring{v}_{a'} \hat{v}_{a'}^{-1})^{\mathcal{B}} )$
 is an automorphism of the induced subgraph on $D_a' \cup D_{a'}'$
 if and only $(id^{\mathcal{B}},(\mathring{v}_{a}^{-1} \mathring{v}_{a'}
 \hat{v}_{a'}^{-1} \hat{v}_a)^{\mathcal{B}})$ is. Since both $\mathring{Z}_p^3$
 and $\hat{Z}_p^3$ are abelian we have
 \[ \left( id^{\mathcal{B}},(\mathring{v}_{a}^{-1} \mathring{v}_{a'}
 \hat{v}_{a'}^{-1} \hat{v}_a)^{\mathcal{B}} \right) =
 \left(id^{\mathcal{B}},(\mathring{v}_{a'} \mathring{v}_{a}^{-1})^{\mathcal{B}}
 (\hat{v}_a \hat{v}_{a'}^{-1})^{\mathcal{B}} \right) \mbox{.} \]
 Clearly, $(\hat{v}_a \hat{v}_{a'}^{-1})^{\mathcal{B} }(a',b_{a'},0)
 =(a,b_{a},0)$ and $(\mathring{v}_{a'} \mathring{v}_{a}^{-1})^{\mathcal{B}}
 (a,b_{a},0)=(a',b_{a'},0)$. Using Corollary \ref{cor1} we get that
 \[ \left(id^{\mathcal{B} },(\mathring{v}_{a'} \mathring{v}_{a}^{-1})^{\mathcal{B}}
 (\hat{v}_a \hat{v}_{a'}^{-1})^{\mathcal{B}} \right)
 =\left(id^{\mathcal{B}},id^{\mathcal{B}} \right) \]
  which is clearly an automorphism on $D_a' \cup D_{a'}'$. This proves that
  $\alpha_3^{\mathcal{B} } \in Aut(\Gamma_1)$.

  If $B_i \sim B_j$, then $\alpha_3 (B_i) \sim \alpha_3(B_j)$ since
   $\alpha_3^{\mathcal{B} } \in Aut(\Gamma_1)$ thus $p_i \in B_i$ is
   connected to $p_j \in B_j$ if and only if $\alpha_3(p_i)$
   is connected to $\alpha_3(p_j)$.

   If $B_i \nsim B_j$, then there exists $a \in \mathbb{Z}_p$
   such that $B_i$ and  $B_j \subset D_a$. Since $\alpha_3$ is
   defined on $D_a$ by an automorphism of $\Gamma$ we have that $p_i \in B_i$ is
   connected to $p_j \in B_j$ if and only if $\alpha_3(p_i)$
   is connected to $\alpha_3(p_j)$, finishing the proof of Lemma \ref{lem9}.
\qed
\end{pf}

Finally, one can see using Lemma \ref{szorzas} \ref{szorzasb} that $\mathring{g}^{\alpha_3} = \hat{g}$.

\subsection{Step 3}\label{step3}

Let us assume that for the generators of the cyclic groups
$\hat{g} \in \hat{\mathbb{Z}}_q$ and $\mathring{g} \in \mathring{\mathbb{Z}}_q$ we have $\mathring{g}= \hat{g}$.

Since $\mathring{g}=\hat{g}$ we have that $\hat{\mathbb{Z}}_p^3$ and
$\mathring{\mathbb{Z}}_p^3$ are
contained in $C_A(\hat{g} )$. Using Sylow's theorem again we may assume
 that $\hat{\mathbb{Z}}_p^3$ and
$\mathring{\mathbb{Z}}_p^3$ are in the same Sylow $p$-subgroup of $C_A( \hat{g} )$.
Using all these assumptions we prove the following Lemma.
\begin{lem}\label{osszefuggo}
\begin{enumerate}
\item\label{het} $\mathring{\mathbb{Z}}_p^3 \times \mathring{\mathbb{Z}}_q \le
\hat{\mathbb{Z}}_q \wr Sym(p^3)$.
\item\label{e} If $\mathring{\mathbb{Z}}_p^3 \times \mathring{\mathbb{Z}}_q \le
\hat{\mathbb{Z}}_q \wr Sym(p^3)$, then for every $\mathring{u}
\in \mathring{\mathbb{Z}}_p^3$ we have $(\mathring{u} ) _b =id$.
\end{enumerate}
\end{lem}
\begin{pf}
\begin{enumerate}
\item $\mathring{\mathbb{Z}}_p^3 \times \mathring{\mathbb{Z}}_q \le
\hat{\mathbb{Z}}_q \wr Sym(p^3)$
since the elements of $\mathring{\mathbb{Z}}_p^3$ and $\hat{g}$ commute.
\item
Let $A' = A \cap \hat{\mathbb{Z}}_q \wr Sym(p^3)$. We have already
 assumed that $\mathring{\mathbb{Z}}_p^3$ and $\hat{\mathbb{Z}}_p^3$
 lie in the same Sylow $p$-subgroup of $A'$, which is generated by 
  $p^3$ disjoint $q$-cycles. Let $\mathring {u}$
be an arbitrary element of $\mathring{\mathbb{Z}}_p^3$.
For every $(b,s) \in \mathbb{Z}_p^3 \times \mathbb{Z}_q$ we
have $\mathring{u}(b,s) = (c, s+t)$ for some $c \in \mathbb{Z}_p^3$
 and $t \in \mathbb{Z}_q$,
where $t$ only depends on $\mathring{u}$ and $b$ since
$\mathring{u} \in \hat{\mathbb{Z}}_q \wr Sym(p^3)$.
The permutation group $\hat{G}$ is transitive, hence there
exist $\hat{u}_1, \hat{u}_2 \in \hat{\mathbb{Z}}_p$ such that
 $\hat{u}_1(0,s)=(b,s)$ and $\hat{u}_2(c,s+t)=(0,s+t)$.
 The order of $\hat{u}_2 \mathring{u} \hat{u}_1$
is a power of $p$ since $\hat{u}_2, \mathring{u}$ and $\hat{u}_1$
 lie in a Sylow $p$-subgroup. Therefore $t =0$ and hence
 $(\mathring{u})_b = id$.
 \end{enumerate}
 \qed
\end{pf}

Lemma \ref{osszefuggo} says that for every $\mathring{u} \in \mathring{\mathbb{Z}}_p^3$ we have $(u)_b =id$.
We use again the graph $\Gamma_1$ defined on $\mathcal{B}$. It is clear that $H_1$ and $H_2$ are regular subgroups in $Aut(\Gamma_1)$
and they are isomorphic to  $\mathbb{Z}_p^3$. Since $\mathbb{Z}_p^3$ is a DCI$^{(2)}$-group \cite{alspachnowitz} we have that there exists
$\mu \in \langle H_1,H_2 \rangle^{(2)}$ such that $H_2 ^{\mu} =H_1$.

Let $\eta =\mu id_{\mathcal{B} }$ be an element of the wreath product $\mathbb{Z}_q \wr Sym(p^3)$. Clearly,
$\eta \in \langle \hat{G},\mathring{ \mathbb{G} } \rangle^{(2)}$ and hence $\eta$ is an
automorphism of $\Gamma_0$, which conjugates $\mathring{\mathbb{Z}}_p^3$ to $\hat{\mathbb{Z}}_p^3$.
Moreover, the base group part of $\eta$ is the identity so $\eta \in C_A(\hat{g})$.
This proves that $\mathring{G}^{\eta} = \hat{G}$, finishing the proof of Theorem \ref{fotetel}.

\end{enumerate}

\end{document}